\documentclass[12pt]{amsart}
\usepackage{amscd}
\usepackage[latin1]{inputenc}
\usepackage{amscd}
\usepackage{latexsym}
\usepackage{pstcol,pst-plot,pst-3d}
\usepackage{multicol}
\psset{unit=0.7cm,linewidth=0.8pt,arrowsize=2.5pt 4}

\newpsstyle{fatline}{linewidth=1.5pt}
\newpsstyle{fyp}{fillstyle=solid,fillcolor=verylight}
\definecolor{verylight}{gray}{0.97}
\definecolor{light}{gray}{0.9}
\definecolor{medium}{gray}{0.85}

\def\NZQ{\Bbb}
\def\NN{{\NZQ N}}

\def\RR{{\NZQ R}}

\let\sect=\cap

\def\F{{\mathcal F}}
\def\D{{\Delta}}
\def\G{{\Gamma}}
\def\a{{\bold a}}
\def\b{{\bold b}}
\def\c{{\bold c}}
\def\e{{\bold e}}

\def\1{{\mathbf 1}}
\def\0{{\mathbf 0}}

\def\opn#1#2{\def#1{\operatorname{#2}}} 
\opn\relint{relint}
\opn\height{ht}
\opn\tr{tr}
\opn\ini{in}
\newtheorem{Theorem}{Theorem}[section]
\newtheorem{Lemma}[Theorem]{Lemma}
\newtheorem{Corollary}[Theorem]{Corollary}
\newtheorem{Proposition}[Theorem]{Proposition}

\newtheorem{Example}[Theorem]{Example}

\begin{document}

\title{Standard graded vertex cover algebras,\\  cycles and leaves}

\author{J\"urgen Herzog, Takayuki Hibi,  Ng\^o Vi\^et Trung and Xinxian Zheng}
\address{J\"urgen Herzog, Fachbereich Mathematik und
Informatik, Universit\"at Duisburg-Essen, Campus Essen, 45117
Essen, Germany} \email{juergen.herzog@uni-essen.de}

\address{Takayuki Hibi, Department of Pure and Applied Mathematics,
Graduate School of Information Science and Technology, Osaka
University, Toyonaka, Osaka 560-0043, Japan}
\email{hibi@math.sci.osaka-u.ac.jp}

\address{Ng\^o Vi\^et Trung, Institute of Mathematics,
Vien Toan Hoc, 18 Hoang Quoc Viet, 10307 Hanoi, Vietnam}
\email{nvtrung@math.ac.vn}

\address{Xinxian Zheng, Fachbereich Mathematik und
Informatik, Universit\"at Duisburg-Essen, Campus Essen, 45117
Essen, Germany} \email{xinxian.zheng@uni-essen.de}

\begin{abstract}
The aim of this paper is to characterize simplicial complexes
which have standard graded vertex cover algebras. This property
has several nice consequences for the squarefree monomial ideals
defining these algebras.
It turns out that such simplicial complexes are closely related to a range of hypergraphs which generalize bipartite graphs and trees. 
These relationships allow us to obtain very general results on standard graded vertex cover algebras which cover  previous major 
results on Rees algebras of squarefree monomial  ideals.
\end{abstract}

\maketitle

\section*{Introduction}

Let  $\Delta$ be a simplicial complex on the vertex set $[n] =
\{1,...,n\}$. Let ${\mathcal F}(\Delta)$ denote the set of the
facets of $\Delta$. An integer vector $\c = (c_1, \ldots, c_n) \in
\NN^n$ is called a {\em cover of order $k$} or a {\em $k$-cover}
of $\Delta$ if $\sum_{i \in F} c_i \geq k$ for all facets $F$ of
$\Delta$. If $\c$ happens to be a $(0,1)$-vector, then $\c$ may be
identified with the subset $C=\{i\in [n] :\; c_i \neq 0\}$ of
$[n]$. It is clear that $\c$ is a  1-cover  if and only if $C$ is
a {\it vertex cover} of $\Delta$ in the classical sense, that is,
$C\sect F\neq \emptyset$ for all $F\in {\mathcal F}(\Delta)$.

Let $S = K[x_1,\ldots,x_n]$ be a polynomial ring in $n$ variables over a field $K$. Let $A_k(\Delta)$ denote the $K$-vector space 
generated by all monomials $x_1^{c_1}\cdots x_n^{c_n}t^k$ such that $(c_1,\ldots, c_n)\in \NN^n$ is a $k$-cover of $\Delta$, where $t$ 
is a new variable. Then
\[
A(\Delta) :=  \bigoplus_{k\geq 0}A_k(\Delta),
\]
is a graded $S$-algebra. We call $A(\Delta)$ the
{\em vertex cover algebra} of $\Delta$ \cite{He}.

Vertex cover algebras have an interesting algebraic interpretation.
For any subset $F \subset [n]$ let $P_F$ denote the ideal of $S$
generated by the variables $x_i$ with $i\in F$. Set
$$I^*(\Delta)  :=\bigcap_{F\in {\mathcal F}(\Delta)}P_F.$$
Then $I^*(\Delta)$ is the ideal generated by the squarefree monomials $x_{i_1}\cdots x_{i_j}$ such that
$\{i_1,\ldots, i_j\}$ is a vertex cover of $\Delta$ and $A(\Delta)$ is
the symbolic Rees algebra of $I^*(\Delta)$. It is shown in \cite{He} that $A(\Delta)$  is a finitely generated, graded and normal 
Cohen-Macaulay $S$-algebra.

Ideals of the form $I^*(\D)$ are exactly the squarefree monomial ideals. A similar notion is the {\it facet ideal} $I(\D)$ which is 
generated by the squarefree monomials $x_{i_1}\cdots x_{i_j}$ such that
$\{i_1,\ldots, i_j\}$ is a facet of $\D$. If we denote by $\D^*$ the simplicial complex whose facets are the minimal vertex covers of 
$\D$, then $I^*(\D^*) = I(\D)$. Thus, one can use vertex cover algebras as a tool for the study of facet ideals.

It is of great interest to know when $A(\Delta)$ is a standard graded algebra, that is, when $A(\Delta)$  is generated over $S$ by 
homogeneous elements  of degree 1. Combinatorially, that means  every $k$-covers  of $\D$ can be written as a sum of $k$ 1-covers for 
all $k \ge 1$. Ideal-theoretically, that means the symbolic powers of $I^*(\D)$ coincide
with the ordinary powers. In particular, this condition implies that the Rees algebra of $I(\D^*)$ is a Cohen-Macaulay normal domain 
and that the associated graded ring is Cohen-Macaulay and reduced.

In case $\Delta$ is pure and 1-dimensional, we identify $\Delta$
with the graph whose edges are the facets of $\Delta$. For such a
simplicial complex it is shown in \cite{He} that $A(\Delta)$ is a
standard graded algebra if and only if $\Delta$ is a bipartite
graph. This result has led us to study more generally the
relationship between standard graded vertex cover algebras and the
combinatorics of the corresponding simplicial complexes.

In the more general situation it proves to
be favorable to identify a given simplicial complex $\Delta$  with
the hypergraph whose edges are the facets of $\Delta$.
Recall that a hypergraph $H$ is a collection of subsets of a finite set of vertices.
The elements of $H$ are called the edges of $H$. Of course, the hypergraphs of
the form ${\mathcal F}(\Delta)$ are special, because there are no
inclusions between the sets representing its edges.
Up to the order of the vertices and the order of the edges, a
hypergraph is determined by its incidence matrix and vice versa.
This matrix is defined as follows: let $x_1,\ldots, x_n$ be the
vertices and $C_1,\ldots, C_m$ be the edges of the hypergraph $H$.
Then the incidence matrix $M=(e_{ij})$ is an $n\times m$ matrix
with  $e_{ij}=1$ if $x_i\in C_j$ and $e_{ij}=0$ if $x_i\not \in
C_j$.

We say that $\D$ is a Mengerian simplicial complex if the incidence matrix of $\F(\D)$ satisfies a certain min-max equation, which is 
known as the Mengerian property in hypergraph theory or as the max-flow min-cut property in integer programming. As a main result in 
Section 1 we show in Theorem \ref{Mengerian 1} that $A(\D)$ is a standard graded algebra if and only if $\D^*$ is a Mengerian 
simplicial complex. The proof of this result is based on the
observation that $A(\Delta)$ is standard graded if and only if
the vertex set of every polarized simplicial complex of $\Delta$ can be decomposed into $k$ vertex covers, where $k$ is the minimum 
number of vertices of the facets. As a consequence, $\D$ is a Mengerian simplicial complex if and only if the symbolic powers of the 
facet ideal $I(\D)$ coincides with the ordinary powers. This result recovers two recent results of Escobar, Villarreal and Yoshino 
\cite{Es} and Gitler, Valencia and Villareal \cite{Gi} on the normality of the Rees algebra of $I(\D)$.

It suggests itself  to ask whether the fact mentioned above that
each bipartite graph has a standard graded vertex cover algebra
may be extended to higher dimensions. This is indeed the case and
is discussed in Section 2. But what is the analogue of a bipartite
graphs in higher dimensions? Recall that a bipartite graph is
characterized by the property that it has no odd cycles. It is
natural to call an alternating sequence of distinct vertices and edges
$v_1,F_1,v_2,F_1,\ldots, v_s,F_s,v_{s+1} = v_1$, $s \ge 2$, in a hypergraph a {\it cycle}
if $v_i,v_{i+1}\in F_i$ for all $i$. Such a
cycle is called {\it special} if no edge contains more the two
vertices of the cycle.

We show in Theorem \ref{balanced}
that for a given simplicial complex $\Delta$, the vertex cover
algebra $A(\Gamma)$ is standard graded for all subcomplexes
$\Gamma\subseteq \Delta$, if and only if $\Delta$ has no special
odd cycles. Hypergraphs having no special odd cycle are called totally
balanced and our proof uses the well-known result of Berge \cite{Be1}
that the vertices of a balanced hypergraph $H$ can be colored by $\min\{|F|:\ F \in H\}$ colors such that every color occurs in every 
edge  of $H$.

Fulkerson, Hoffman and Oppenheim \cite{Fu} showed that a balanced
hypergraph is Mengerian. This implies at once that $A(\Delta^*)$
is standard graded if $\Delta$ has no special odd cycle. This, in
turn,  immediately yields the result of Simis, Vasconcelos and
Villarreal \cite{Si} according to which the edge ideal of a bipartite
graph is normally torsionfree.

A well-known class of balanced hypergraphs are the unimodular
hypergraphs which are defined by the property that all
determinants of the incidence matrix equal $0,\pm 1$. We say that
$\Delta$ is unimodular if ${\mathcal F}(\Delta)$ unimodular. In
this case, the above mentioned results  imply that $A(\Delta)$ and
$A(\Delta^*)$ are standard graded algebras.

A facet $F$ of a simplicial complex of $\D$ is called a
{\em leaf} if either $F$ is the only facet of $\Delta$, or there
exists  $G\in {\mathcal F}(\Delta)$, $G\neq F$ such that $H\sect
F\subset G\sect F$ for each $H\in{\mathcal F}(\Delta)$ with $H\neq
F$.  If each subcomplex $\Gamma$ of $\D$ has a leaf, then $\Delta$
is called a forest. If, in addition $\D$ is connected, then $\D$
is called a tree. These notions were introduced by Faridi \cite{Fa} who proved, among other things, that
the Rees algebra of the facet ideal of a tree is a Cohen-Macaulay domain.

In Theorem \ref{forest} we show that $\Delta$ is a forest if and only
if $\Delta$ has no special cycle of length $\geq 3$.
Hypergraphs having no special cycle of length $\geq 3$ are called totally
balanced and our proof is based on the characterization of totally balanced hypergraphs in terms of the so-called greedy matrix found 
by Hoffman, Kolen and Sakarovitch \cite{Ho} and Lubiw \cite{Lu}.
The results in Section 2 imply that $A(\Delta)$ and $A(\Delta^*)$
are standard graded if $\Delta$ is a forest. This then implies the
above mentioned result of Faridi and stronger assertions on the
facet ideal of a tree.

There is another interesting property  of forests, shown in
Corollary \ref{good leaf 1}, namely that each forest has a good
leaf, that is, a facet which is a leaf of each subcomplex to which
it belongs. The notion good leaf was introduced in the thesis of
Zheng \cite{Zh}. The existence of good leaves implies immediately
that each forest has a good leaf order, that is to say, that the
facets $F_1,\ldots, F_m$ of $\D$ can be ordered in such a way that
$F_i$ is a good leaf of the subcomplex whose facets are
$F_1,\ldots, F_i$. It was recently shown by  Pelsmajer, Tokaz and
 West \cite{Pe} that totally balanced hypergraphs
have good leaf orders. Their proof compared with the arguments
given here is relatively complicated.

Using good leaf orders and  an algebraic result of Conca and De
Negri \cite{Co} we show in Corollary \ref{Rees CM} that
if $I$ is a graded ideal in a polynomial ring, whose initial ideal
with respect to a suitable term order is the facet ideal of a
forest, has the property that the Rees algebra of $I$ is
Cohen-Macaulay and that the associated graded ring of $I$ is
reduced.

We close Section 3 by showing in Theorem \ref{good leaf 2} that if
$\Delta$ is a simplicial complex, $F$ is a good leaf of $\Delta$
and $\Gamma$ is the simplicial complex obtained from $\Delta$ by
removing $F$, then the highest degree of the generators of
$A(\Delta)$ is the same as the highest degree of  the generators of
$A(\Gamma)$. We use this result in Section 4 to study
vertex cover algebras of quasi-forests.

A quasi-forest is a simplicial complex whose facets can be ordered
$F_1,\ldots, F_m$ such that for all $i$, $F_i$ is a leaf of the
simplicial complex with the facets $F_1,\ldots, F_i$.  It is clear
that any forest is a quasi-forest. The significance of
quasi-forests results from a theorem of Dirac  \cite{D} who proved
that (translated into the language of quasi-forests)  a simplicial
complex is a quasi-forest if and only if its 1-skeleton is a
chordal graph.

There is a subclass of the chordal graphs, called strongly chordal. Farber
\cite{Farber} showed that strongly chordal graphs can be described
by forbidden induced subgraphs, the so-called trampolines, and that
a simplicial complex is a forest if and only if its 1-skeleton is a strongly chordal graphs.
Comparing this result with Dirac's theorem the difference between forest and
quasi-forests becomes apparent and one would expect that unlike
forests the quasi-forests have non-standard graded vertex cover
algebras. However, this is not always the case as can be shown by
examples. On the other hand, we show in our final  Theorem \ref{also}
that a quasi-forest with the property that it is connected in
codimension 1 and that each face of codimension 1 belongs to at most
two facets is a forest if and only if
$A(\Delta)$ is standard graded.

We have seen in this paper several applications of hypergraph theory on vertex cover algebras and facet ideals.
It would be of interest to use algebraic methods to solve problems in hypergraph theory.
Moreover, the notion of vertex cover algebras has been introduced for weighted simplicial complexes.
It remains to see whether one can extend some of the results of this paper for general vertex cover algebras.

\section{Standard vertex cover algebras}

We adhere to the notions of the introduction.

Let  $\Delta$ be a simplicial complex on the vertex set $[n]$. We
say that a $k$-cover $\c$ of $\D$ is {\em decomposable} if there
exists an $i$-cover $\a$ and a $j$-cover $\b$ such that $\c=
\a+\b$, $k=i+j$ and $a, b\neq 0$. If $\c$ is not decomposable, we
call it {\em indecomposable}. It is clear that the indecomposable
covers of order $\ge 1$ correspond to a minimal homogeneous set of
generators of the $S$-algebra $A(\Delta)$.

An obvious vertex cover of  $\D$ is the {\it full cover} $[n]$, which corresponds to the vector $(1,...,1)$.
The maximal order of this cover is obviously given by the number
$$s(\D) = \min\{|F|:\ F \in \F(\D)\}.$$
If $A(\D)$ is a standard graded $S$-algebra, then $[n]$ can be decomposed into $s(\D)$ vertex covers of $\D$. In this case, $\D$ is 
called {\it totally decomposable}. We shall use this property  to give a characterization of standard graded vertex cover algebras.

Let $\c = (c_1,...,c_n) \in \NN^n$ be an arbitrary integral vector. We associate with $\c$ a new set of vertices
$$X^\c = \{x_{ij}:\ i = 1,...,n, j = 1,...,c_i\}.$$
To each subset $F$ of $[n]$ we associate the subset
$$F^\c := \{x_{ij}:\ i \in F, j = 1,...,c_i\}$$
of $X^c$. Let $\Delta^\c$ denote the simplicial complex on $X^\c$ whose facets are the minimal sets of the form  $F^\c$, $F \in 
\F(\Delta)$. Following the method of polarization of monomials (see e.g. \cite{Br}) we call $\Delta^\c$ the {\it polarization} of 
$\Delta$ with respect to $\c$.

The following result allows us to reduce the decomposition of  arbitrary covers of $\D$ to the decomposition of the full cover of their  
polarizations.

\begin{Lemma} \label{polar-1}
An integer vector $\c \in \NN^n$ can be written as a sum of $k$ 1-covers if and only if $X^\c$ can be decomposed into $k$ vertex covers 
of $\D^\c$.
\end{Lemma}

\begin{proof}
Assume that $\c = \c_1 + \cdots + \c_k$, where $\c_1,...,\c_k$ are 1-covers of $\D$.
For $h = 1,...,k$ let $\c_h = (c_{h1},...,c_{hn})$ and
$$C_h := \{x_{ij}:  i \in 1,...,n,\ j = c_{1i}+ \cdots + c_{h-1i}+1,...,c_{1i}+ \cdots + c_{h-1i}+c_{hi}\}.$$
It is easily seen that $X^\c = C_1 \cup \cdots \cup C_k$ and that $C_1,...,C_k$ are disjoint vertex covers of $\D^\c$.
\par

Conversely, assume that $X^\c = C_1 \cup \cdots \cup C_k$ and that $C_1,...,C_k$ are disjoint vertex covers of $\D^\c$. For $h = 
1,...,k$ let $\c_h = (c_{h1},...,c_{hn})$ with $c_{hi} := |x_{ij}:\ x_{ij} \in C_h|$, $i =1,...,n$.
It is obvious that $\c = \c_1 + \cdots + \c_k$. Each $\c_h$ is a 1-cover of $\D$ because for all facets $F$ of $\D$,
$F^\c$ contains a facet of $\D^\c$. But $C_h$ contains a vertex, say $x_{ij}$, of this facet of $\D^\c$. Therefore,
$i \in F$ and $c_{hi} \ge 1$, which implies $\sum_{i \in F}c_{hi}\ge 1$.
\end{proof}

\begin{Proposition} \label{full cover}
The vertex cover algebra $A(\D)$ is standard graded if and only if $\D^\c$ is totally decomposable for all $\c \in \NN^n$.
\end{Proposition}

\begin{proof}
$A(\D)$ is standard graded if and only if every $k$-covers $\c$ of $\D$ can be written as a sum of $k$ 1-covers for all $k \ge 1$. By 
Lemma \ref{polar-1}, this is equivalent to the condition that $X^\c$ can be decomposed into $k$ vertex covers of $\D^\c$. It is easily 
seen that
$$s(\D^\c) = \min\{|F^\c|:\ F \in \F(\D)\} = \min\{\sum_{i \in F}c_i:\ F \in \F(\D)\}$$
is the maximal order of $\c$. Therefore, $A(\D)$ is standard graded if and only if $X^\c$ can be decomposed into $s(\D^\c)$ vertex 
covers of $\D^\c$.
\end{proof}

Let $o(\c)$ denote the maximal order of $\c$ and $\sigma(\c)$ the maximum number $k$ such that $\c$ can be written as a sum of $k$ 
1-covers of $\D$. Then $o(\c) = s(\D^\c)$ and $\sigma(\c)$ is  the maximum number $k$ such that $X^\c$ can be decomposed into $k$ 
vertex covers of $\D^\c$. Thus, $\sigma(\c) \le o(\c)$ and $\D^\c$ is totally decomposable means $\sigma(\c) = o(\c)$.

The above notions have the following meanings in hypergraph theory.
Recall that a hypergraph $H$ is a collection of subsets of a vertex set. The elements of $H$ is called the edges of $H$. One calls a 
set of vertices  meeting all edge of $H$ a {\it transversal} of $H$.
If $H$ has a partition into $k$ transversals, where $k = \min\{|F|:\ F \in H\}$, then $H$ is said to have the Gupta property 
\cite[3.1]{Du}.  Therefore, $\D^\c$ is totally decomposable if and only if  the hypergraph $\F(\D^\c)$ has the Gupta property.

Berge  \cite[Chapter 5, Lemma, p.~207]{Be2} already studied hypergraphs $H$ such that every polarization of $H$ has the Gupta property, 
using different notations.  Following his approach we can give a characterization of standard graded algebras by means of a min-max  
property.

Let $C_1,...,C_m \subset X$ be the minimal vertex covers of $\D$.
Let $\D^*$ denote the simplicial complexes whose facets are $C_1,...,C_m$. It is well-known that
$$(\D^*)^* = \D.$$
Let $M$ be the facet-vertex incidence matrix of $\D^*$.
Then we have the following formulas for $o(\c)$ and $\sigma(\c)$ in terms of $M$.

\begin{Lemma} \label{incidence}
Let $\1 $ denote the vector $(1,...,1)$ of  $\NN^m$. Then
\par {\rm (i) }  $o(\c) = \min\{\a \cdot \c:\  \a \in \NN^n,\  M\cdot\a \ge \1 \}$,
\par {\rm (ii)}  $\sigma(\c) = \max\{\b \cdot \1 :\ \b \in \NN^m,\ M^T\cdot\b \le \c\}$.
\end{Lemma}

\begin{proof}
The vectors $\a \in \NN^n$ with  $M\cdot\a \ge \1 $ are the 1-covers of $\D^*$.
Since the minimal 1-covers of $\D^*$ correspond to the facets of $\D$, we have
\begin{align*}
o(\c) & = \min\{\sum_{i \in F}c_i:\  F \in \F(\D)\} \\
 & = \min\{\a \cdot \c:\  \a \in \{0,1\}^n,\  M\cdot\a \ge \1 \}\\
& = \min\{\a \cdot \c:\  \a \in \NN^n,\  M\cdot\a \ge \1 \}.
\end{align*}

Let $\c_1,...,\c_m$ be the (0,1) vectors of  $C_1,...,C_m$. As $\c_1,...,\c_m$ are the minimum 1-covers of $\D$,
$\c$ can be written as a sum of $k$ 1-covers of $\D$ if there exist $b_1$ copies of $\c_1$,..., $b_m$ copies of $\c_m$ such that $k = 
b_1 + \cdots + b_m$ and $b_1\c_1 + \cdots + b_m\c_m \le \c$. Let $\b = (b_1,...,b_m)$.
Then $\b\cdot\1  = b_1 + \cdots + b_m$ and
$M^T\cdot\b = b_1\c_1 + \cdots + b_m\c_m$. Therefore,
\begin{align*}
\sigma(\c) & = \max\{k:\ \text{$\c$ can be written as a sum of $k$ 1-covers of $\D$}\}\\
& = \max\{\b\cdot\1 :\  \b \in \NN^m,\ M^T\cdot\b \le \c\}.
\end{align*}
\end{proof}

Following the terminology of   hypergraph theory \cite[Chapter 7]{Be2}, we call $\Delta^*$ a {\it Mengerian
simplicial  complex} if for all $\c \in \NN^n$,
$$\min\{\a \cdot \c:\  \a \in \NN^n,\  M\cdot\a \ge \1 \}
= \max\{\b \cdot \1 :\ \b \in \NN^m,\ M^T\cdot\b \le \c\}.$$
This definition applies to an arbitrary simplicial complex or a hypergraph with suitably adopted notation.
The Mengerian property is also known under the name max-flow min-cut property in integer linear programming.

Lemma \ref{incidence} immediately yields the following characterization of standard graded vertex cover algebras.

\begin{Theorem}\label{Mengerian 1}

The vertex cover algebra $A(\D)$ is standard graded  if and only if $\D^*$ is a Mengerian simplicial complex.
\end{Theorem}

As $\D = (\D^*)^*$, this result can be reformulated as follows.

\begin{Corollary} \label{cover complex}
$\D$ is a Mengerian simplicial complex if and only if $A(\D^*)$ is
a standard graded algebra.
\end{Corollary}

Let $I(\D)$ be the {\it facet ideal} of $\D$, the ideal generated by the monomials of the variables of the facets of $\D$. Then $I(\D) 
= I^*(\D^*)$. This interpretation of $I(\D)$ leads to interesting relationships between the Mengerian property and properties of the 
facet ideal.

\begin{Corollary} \label{facet ideal 1}
The following conditions are equivalent:
\par {\rm (i) }  $\D$ is a Mengerian simplicial complexes.
\par {\rm (ii) } $I(\D)^{(k)} = I(\D)^k$ for all $k \ge 0$,
\par {\rm (iii)} The associated graded ring of $I(\D)$ is reduced,
\par {\rm (iv)} $I(\D)$ is normally torsionfree (i.e. all powers of $I(\D)$ have the same associated prime ideals).
\par \noindent
 Moreover, if one of these conditions is satisfied, then the Rees algebra of $I(\D)$ is a normal Cohen-Macaulay domain and the 
associated graded ring is Cohen-Macaulay.
\end{Corollary}

\begin{proof}
The equivalence of (i) and (ii) follows from Corollary \ref{cover complex} because  $A(\D^*)$ is the symbolic Rees algebra of $I(\D)$. 
The equivalence of (ii) to (iii) and (iv) is well-known \cite{HuV}.
These conditions imply that the Rees algebra of $I(\D)$ coincides with $A(\D^*)$.
By \cite[Theorem 4.3]{He}, $A(\D^*)$ is a normal Cohen-Macaulay domain. It is known that the Cohen-Macaulay property of the Rees 
algebra implies the Cohen-Macaulay property of the associated graded ring \cite {Hu}.
\end{proof}

\begin{Corollary}\label{fiber}
Assume that $\D$ is  pure and Mengerian. Let $k[I(\D)]$ denote the toric ring generated by the monomial generators of $I(\D)$. Then 
$k[I(\D)]$ is a normal Cohen-Macaulay domain.
\end{Corollary}

\begin{proof}
Since $\D$ is pure, $I(\D)$ is generated by monomials of the same degree, say $d$. Therefore, we may view  the Rees algebra ${\mathcal 
R}$ of $I(\D)$ as a bigraded $k$-algebra with ${\mathcal R}_{(h,k)} = (I(\D)^k)_{dh}t^k$ for $h, k \ge 0$ and $k[I(\D)]$ as the 
subalgebra $\oplus_{k \ge 0}{\mathcal R}_{(k,k)}$. Since $\D$ is Mengerian, $A(\D)$ is a normal domain by Corollary \ref{facet ideal 
1}. Hence, $k[I(\D)]$ is also a normal domain. By \cite{Hoch}, this implies that $k[I(\D)]$ is Cohen-Macaulay.
\end{proof}

The relationship between facet ideals and the Mengerian simplicial complexes was already studied in a recent paper of Gitler, Valencia 
and  Villareal  \cite[Theorem 3.5]{Gi}. They proved that $\D$ is Mengerian if and only if  the polyhedron
$$Q(\D) := \{\a \in \RR^n:\  \a \ge \0,\ M \cdot \a \ge \1 \}$$
($M$ is now the incidence matrix of $\D$ and $\0$ is the vector of zero components) has integral vertices and the Rees algebra of 
$I(\D)$ is a normal domain. On the other hand, Escobar, Villarreal and Yoshino \cite[Proposition 3.4]{Es} showed that the latter 
conditions are satisfied if and only if  the associated graded ring of $I(\D)$ is reduced. Hence, Corollary \ref{facet ideal 1} can be 
deduced from their results.

It is not hard to see that $Q(\D)$ has integral vertices if and only if the normalization of the Rees algebra of $I(\D)$ coincides with 
the symbolic Rees algebra. Therefore, one can also recover the two afore-mentioned results from Corollary \ref{facet ideal 1}.

\begin{Example}\label{Mengerian 2}
{\rm Let $\D$ be the simplicial complex of Picture \ref{K4} which has the facets
$$\{1,2,3\},\{1,5,6\},\{2,4,6\},\{3,4,5\}.$$

\begin{figure}[hbt]
\begin{center}
\psset{unit=1.5cm}
\begin{pspicture}(-2,-1.1)(2,1.2)
\pspolygon[style=fyp,fillcolor=medium](0,1)(-0.866,0.5)(-0.866,-0.5)
 \pspolygon[style=fyp, fillcolor=medium](0,1)(0.866,0.5)(0.866,-0.5)
 \pspolygon[style=fyp, fillcolor=medium](0,-1)(-0.866,-0.5)(0.866,-0.5)
  \pspolygon[style=fyp, fillcolor=medium](0.866,0.5)(0,-1)(-0.866,0.5)
 \psline(-0.866,-0.5)(0.866,-0.5)
 \psline(-0.866,-0.5)(0,1)
 \psline(0,1)(0.866,-0.5)
 \rput(0,1.2){1}
 \rput(1.1,0.5){2}
 \rput(1.1,-0.5){3}
 \rput(0,-1.2){4}
 \rput(-1.1,-0.5){5}
 \rput(-1.1,0.5){6}
\end{pspicture}
\end{center}
\caption{}\label{K4}
\end{figure}

Then $\D^*$ is the simplicial complex with the facets
$$\{1,4\},\{2,5\},\{3,6\},\{1,2,3\},\{1,5,6\},\{2,4,6\},\{3,4,5\}.$$
It is known that $\D$ is not Mengerian but $\D^*$ is Mengerian \cite[p.~198 and p.~209]{Be2}.
We will use Proposition \ref{Mengerian 1} to establish this fact.\par

$\D$ is not Mengerian because $A(\D^*)$ is not standard graded. In fact, the full cover of $\D^*$ corresponds to  an indecomposable 
2-cover of $\D$, which implies that $A(\D)$ has a generator in degree 2.

$\D^*$ is Mengerian because $A(\D)$ is standard graded. Let $\c = (c_1,...,c_6)$ be an arbitrary $k$-cover of $\D$ of maximal order $k 
\ge 2$. It suffices to show that $\c$ can written as a sum of an 1-cover and a $(k-1)$-cover. Let $f = x_1^{c_1}\cdots x_6^{c_6}$. If 
$f$ is divisible by a monomials $g$ of the forms $x_1x_4, x_2x_5,x_3x_6$, then $f/g$ correspond to a cover of order $k-1$ because $g$ 
correspond to a 1-covers of $\D$ which meets every edge at only one vertex. Therefore, the decomposition $f = (f/g)g$ corresponds to a 
partition of $\c$ into an 1-cover and a $(k-1)$-cover of $\D$. If $f$ is not divisible by the monomials $x_1x_4, x_2x_5,x_3x_6$, then 
the support of $f$ must be one of the sets $\{x_1,x_2,x_3\},\ \{x_1,x_5,x_6\},\ \{x_2,x_4,x_6\},\ \{x_3,x_4,x_5\}.$
Without restriction we may assume that the support of $f$ is $\{x_1,x_2,x_3\}$. Using the fact
$$f \in (x_1,x_2,x_3)^k\cap (x_1,x_5,x_6)^k \cap (x_2,x_4,x_6)^k \cap(x_3x_4x_5)^k$$
we can show that $f$ is divisible by $x_1^kx_2^kx_3^k$. Since $x_1x_2x_3$ and $x_1^{k-1}x_2^{k-1}x_3^{k-1}$ corresponds to a 1-cover 
and a $(k-1)$-cover of $\D$, $\c$ can written as a sum of an 1-cover and a $(k-1)$-cover of $\D$.}
\end{Example}

\section{Cycles of simplicial complexes}

Let $\D$ be a simplicial complex. A subcomplex of $\D$ is a simplicial complex $\G$ such that the facets of $\G$ are facets of $\D$, 
denoted by $\G \subseteq \D$. The aim of this section is to study simplicial complexes for which the vertex cover algebra of every 
subcomplex is standard graded.

We are inspired of the result  of \cite[Theorem 5.1]{He} (see also
\cite[Corollay 2.6]{Gi1}) that the vertex cover algebra of a
simple graph is
standard graded if and only if the graph is bipartite. This implies that the vertex cover algebra of every subgraph of a bipartite 
graph is standard graded because subgraphs of a bipartite graph are also bipartite. It is well-known that a simple graph is bipartite 
if and only if it has no odd cycle. It turns out that the notion of cycle is the right tool to characterize the above kind of 
simplicial complexes.

A {\it cycle} or, more precisely, an $s$-cycle of $\D$ ($s \ge 2$) is an alternating sequence  of distinct vertices and facets 
$v_1,F_1,...,v_s,F_s,v_{s+1}=v_1$  such that  $v_i, v_{i+1} \in F_i$ for $i = 1,...,s$. A cycle is {\it special} if it has no facet 
containing more than two vertices of the cycle.  Notice that a cycle of a graph is always special.

In the simplicial complex of Figure \ref{triangle}, the cycle $1,\{1,2,4\},2,\{2,3,4\},3,\{1,3,4\},1$ is special, whereas the cycle 
$1,\{1,2,4\},2,\{2,3,4\},4,\{1,3,4\},1$ is not.

\begin{figure}[hbt]
\begin{center}
\psset{unit=1cm}
\begin{pspicture}(0,1)(2,3.2)
\pspolygon[style=fyp,fillcolor=medium](0,1)(2,1)(1,1.58)
\pspolygon[style=fyp, fillcolor=medium](0,1)(1,1.58)(1,2.72)
\pspolygon[style=fyp, fillcolor=medium](2,1)(1,1.58)(1,2.72)
\rput(1.3,1.75){4} \rput(1, 3){1} \rput(2.2,0.8){3} \rput(-0.2,0.8){2}
\end{pspicture}
\end{center}
\caption{}\label{triangle}
\end{figure}

Our starting point is the following simple observation.

\begin{Lemma} \label{2-cover}
Assume that $A(\G)$ has no generators in degree 2 for all $\G \subseteq \D$. Then $\D$ has no special odd cycle.
\end{Lemma}

\begin{proof}
Assume to the contrary that $\Delta$ has a special cycle $ v_1,F_1,...,v_s,F_s,v_1$, where $s$ is an odd number.
Let $\Gamma$ be the subcomplex of $\D$ with the facets $F_1,...,F_s$. Then  $C = \{v_1,...,v_s\}$ corresponds to a 2-cover of $\Gamma$. 
By assumption,  there are two disjoint  vertex covers $C_1$ and $C_2$ of $\Gamma$ such that $C = C_1 \cup C_2$.  Since the cycle is 
special, every facet of $\Gamma$ has exactly two vertices in $C$. Therefore, one of the two vertices belongs to $C_1$ and the other to 
$C_2$. It follows that the numbers $|C_1|$ and $|C_2|$ of vertices of $C_1$ and $C_2$ are equal. Hence $s = |C_1| + |C_2|$ is an even 
number, a contradiction.
\end{proof}

The above observation led us to the following result.

\begin{Theorem} \label{balanced}
The following conditions are equivalent:
\par {\rm (i) } The vertex cover algebra $A(\Gamma)$ is standard graded for all $\Gamma \subseteq \Delta$,
\par {\rm (ii) } The vertex cover algebra $A(\Gamma)$ has no generator of degree 2  for all $\Gamma \subseteq \Delta$,
\par {\rm (iii) } $\Delta$ has no special odd cycle.
 \end{Theorem}

It suffices to show that if $\Delta$ has no special odd cycle, then $A(\Delta)$ is a standard graded algebra. In fact, since the 
assumption implies that all $\Gamma \subseteq \Delta$ have no special odd cycle, it also implies that $A(\Gamma)$ is a standard graded 
algebra.

By Proposition \ref{full cover} we have to show that the polarization $\D^\c$ is totally decomposable for all $\c \in \NN^n$.

\begin{Lemma} \label{polar-2}
If $\Delta$ has no special odd cycle, then so does $\Delta^\c$.
\end{Lemma}

\begin{proof}
It suffices to show that every special cycle of $\D^\c$ of length $\ge 3$ corresponds to a special cycle of $\D$ of the same length. 
Let $z_1,F_1^*,...,z_s,F_s^*,z_1$ be a special cycle of $\Delta^\c$. For $k = 1,...,s$ let $v_k = i$ if $z_k = x_{ij}$. Let $Q$ denote 
the sequence $v_1,F_1,...,v_s,F_s,v_1$.

We first show that $v_1,...,v_s$ are different vertices. Assume to the contrary that $v_h = v_k$ for some indices $h < k \le s$. By the 
defintion of $\D^\c$,  $z_h,z_{k-1},z_k\in F_{k-1}^\c$ and $z_h,z_k,z_{k+1} \in F_k^\c$.
Since $z_h \neq z_k$ and since $z_h$ can coincides with only one of the vertices $z_{k-1},z_{k+1}$, one of the facets $F_{k-1}^\c, 
F_k^\c$ must contain three different vertices of $P$, a contradiction. \par

It is obvious from the definition of $\D^\c$ that $F_1,...,F_s$ are different, $v_i,v_{i+1} \in F_i$ for $i = 1,...,s$, and no facet of 
$Q$ contains three different vertices of $Q$. Therefore, $Q$ is a special $s$-cycle of $\Delta$.
\end{proof}

By the above discussion, Theorem \ref{balanced} now follows from
the following well-known result(\cite{Be1} and \cite[Corollary 2,
p. 177]{Be2})  of Berge  in hypergraph theory.

\begin{Theorem}\label{Berge}
Assume that $\D$ has no special odd cycle. Then $\D$ is totally decomposable.
\end{Theorem}

Notice that a hypergraph $H$ is called {\it balanced} if $H$ has no special odd cycle.
By this definition, $\D$ has no special odd cycle means the hypergraph $\F(\D)$ is balanced.

The notion of balanced simplicial complex has been used for an other property. Stanley  \cite[4.1]{Sta} called $\D$
a {\it balanced simplicial complex} if $\D$ has a coloring of the vertices by $\dim \D + 1$ colors such that the vertices of every 
facet has different colors. According to \cite{Be1} \cite[Corollary 1, p.177]{Be2}, $\D$ has no special odd cycle if and only if every 
subcomplex of $\D$ is balanced in the sense of Stanley. The balanced simplicial complex of Figure \ref{hole} has a special 3-cycle.

\begin{figure}[hbt]
\begin{center}
\psset{unit=1cm}
\begin{pspicture}(0,0.5)(3,2.5)
\pspolygon[style=fyp, fillcolor=medium](0.5,0.5)(1.5,0.5)(1,1.4)
 \pspolygon[style=fyp, fillcolor=medium](1.5,0.5)(2,1.4)(2.5,0.5)
 \pspolygon[style=fyp, fillcolor=medium](1,1.4)(1.5,2.3)(2,1.4)
\end{pspicture}
\end{center}
\caption{}\label{hole}
\end{figure}
\medskip

Balanced hypergraphs have several interesting characterizations. For instance, a hypergraph is balanced if and only if every partial 
subgraph satisfies the K\"onig property, that is, the minimum number of vertices in a transversal is equal the maximum number of 
disjoint edges (Berge and La Vergnas \cite{Be3}).

In particular, Fulkerson, Hoffman and Openheim \cite{Fu} showed that balanced hypergraphs are Mengerian. As a consequence, simplicial 
complexes without special odd cycle are Mengerian. By Corollary \ref{cover complex}, this implies the following result on the 
simplicial complex $\D^*$ of the minimal vertex covers of $\D$.

\begin{Theorem}\label{facet ideal 2}
Assume that $\D$ has no special odd cycle. Then $A(\D^*)$ is a standard graded algebra.
\end{Theorem}

Notice that by Corollary \ref{facet ideal 1}, there are other interesting consequences on the facet ideal $I(\D)$. As bipartite graphs 
are exactly graphs without odd cycle, we immediately obtain the
following result \cite{Si} of Simis, Vasconcelos and Villarreal.

\begin{Corollary} \label{Si}
Let $G$ be a bipartite graph. Then the edge ideal $I(G)$ is normally torsionfree.
\end{Corollary}

\begin{Example}
{\rm Let $\D$ be the simplicial complex with the facets
$$\{1,2,3\},\{1,5,6\},\{2,4,6\},\{3,4,5\}.$$
We have seen in Example \ref{Mengerian 2} that $A(\D)$ is standard graded. But $A(\G)$ is not standard graded for the subcomplex $\G 
\subseteq \D$ with the facets   $\{1,2,3\},\{1,5,6\},\{3,4,5\}$ which has the special odd cycle
$1,\{1,2,3\},3,\{3,4,5\},5,\{1,5,6\},1$. This cycle is also a special odd cycle of the simplicial complex $\D^*$ which has the facets
$$\{1,4\},\{2,5\},\{3,6\},\{1,2,3\},\{1,5,6\},\{2,4,6\},\{3,4,5\}.$$
Since $\D^*$ is Mengerian, this shows that the class of simplicial complexes without special odd cycle is properly contained in the 
class of Mengerian simplicial complexes.}
\end{Example}

A simplicial complexes without special odd cycle can be also characterized in terms of its incidence matrix. In fact, a special cycle 
corresponds to an $s \times s$ submatrix of the form
$$\begin{pmatrix}
1 & 0 & 0 & \cdot & \cdot  & 0 & 1\\
1 & 1 & 0 & \cdot & \cdot & 0 & 0\\
0 & 1 & 1 & \cdot & \cdot & 0 & 0\\
\cdot &  & & & & &\cdot\\
 \cdot &  & & & &&  \cdot\\
\cdot &  & & & & 1 & 0\\
0 &\cdot & \cdot & \cdot & 0  & 1 & 1
\end{pmatrix}$$
with $s \ge 2$. Therefore, $\D$ has no special odd cycle if and only if its incidence matrix has no such $s \times s$ submatrix with 
odd $s$ even after a permutation of rows and columns.

We say that $\D$ is an {\it unimodular} simplicial complex if every square submatrix of its incidence matrix has determinant equal to 
$0,\pm 1$.

The above matrix has determinant equal to $2$ if $s$ is odd. Therefore, an unimodular simplicial complex has no special odd cycle. By 
Theorem \ref{balanced} and Theorem \ref{facet ideal 2} we obtain the following consequence.

\begin{Corollary}
Let $\D$ be an unimodular simplicial complex. Then $A(\D)$ and $A(\D^*)$ are standard graded algebras.
\end{Corollary}

Notice that simplicial complexes without odd cycle are unimodular \cite[Theorem 5, p.~164]{Be2} and that a simplicial complex of 
dimension $\le 2$ is unimodular if and only if it has no special odd cycle \cite[Corollary, p.~176]{Be2}. Therefore, bipartite graphs 
are exactly unimodular graphs.

\section{Leaves of simplicial complexes}

Let $\D$ be a simplicial complex. A facet $F$ of $\D$ is called a {\em leaf}, if
either $F$ is the only facet of $\Delta$, or there exists  $G\in {\mathcal F}(\Delta)$, $G\neq F$ such that $H\sect F\subset G\sect F$ 
for each $H\in{\mathcal F}(\Delta)$ with $H\neq F$.

The simplicial complex in Figure \ref{Fig2} has two leaves, namely $\{1,2,3\}$ and $\{4,5\}$.

\begin{figure}[hbt]
\begin{center}
\psset{unit=1.5cm}
\begin{pspicture}(3,0.8)(7,2.3)
\pspolygon[style=fyp,fillcolor=medium](3.8,1.5)(4.7,2)(4.7,1)
 \pspolygon[style=fyp, fillcolor=medium](4.7,1)(4.7,2)(5.6,1.5)
 \psline(5.6,1.5)(6.6,1.5)
\rput(3.6,1.5){1} \rput(4.7, 2.2){2} \rput(4.7,0.8){3}
\rput(5.7,1.3){4} \rput(6.8,1.5){5} 
\end{pspicture}
\end{center}
\caption{}\label{Fig2}
\end{figure}

A simplicial complex $\D$ is called a {\em forest}, if each subcomplex $\Gamma$ of $\Delta$  has a leaf. A forest is called a {\em 
tree} if it is connected \cite{Fa}.

The simplicial complex in Figure \ref{Fig2} is a tree. However
the simplicial complex in Figure \ref{triangle} is not a tree.

We will characterize forests by means of the notion of special cycle introduced in the preceding section. For that we need the 
following notation.

Let $\a$ and $\b$ be two vectors of integers of the same length. We define $\a \prec \b$ if the right-most nonzero component of $\a - 
\b$ is negative. Let $M$ be a matrix of integers with rows $\a_1,...,\a_m$ and column $\b_1,...,\b_n$. We say that $M$ is {\it 
canonical} if  $\a_1 \prec \cdots \prec \a_m$ and $\b_1 \prec \cdots \prec \b_n$. A {\em
canonical form} of $M$ is a canonical matrix which can be obtained by permuting the rows and the columns of $M$. The following result 
is well-known in hypergraph theory and linear algebra.

\begin{Lemma}\label{canonical}
Every matrix of integers has a canonical form.
\end{Lemma}

\begin{proof}
Let $M = (a_{ij})$ be a $m \times n$ matrix of integer. Set $\delta_k = \sum_{i+j=k} a_{ij}$ and consider the vector $\delta(M) = 
(\delta_2, \ldots, \delta_{m+n})$.  If $M$ has two rows or columns $\a_{i_1}$ and $\a_{i_2}$, $i_1 < i_2$, in the wrong order $\a_{i_2} 
\prec \a_{i_1}$ and if $M'$ is the matrix obtained by permuting $\a_{i_1}$ and $\a_{i_2}$, then $\delta(M) \prec \delta(M')$. Hence 
permuting the rows and the columns of $M$ to maximize $\delta(M)$ will yield a canonical form of $M$.
\end{proof}

For example, the incidence matrix of a special cycle has the canonical form
$$\begin{pmatrix}
1 & 1 & 0 & \cdot & \cdot  & 0 & 0\\
1 & 0 & 1 & \cdot & \cdot & 0 & 0\\
0 & 1 & 0 & \cdot & \cdot & 0 & 0\\
\cdot &  & & & & &\cdot\\
 \cdot &  & & & &&  \cdot\\
\cdot &  & & & & 0 & 1\\
0 &\cdot & \cdot & \cdot & 0  & 1 & 1
\end{pmatrix}.$$

\begin{Theorem} \label{forest}
The following conditions are equivalent:
\par {\rm (i) } $\D$ is a forest,
\par {\rm (ii) } $\D$ has no special cycle of length $\ge 3$,
\par {\rm (iii) } The incidence matrix of $\D$ has a canonical form which contains no submatrix of the form
\[
B = \begin{pmatrix}
1 & 1\\
1 & 0
\end{pmatrix}.
\]
\end{Theorem}

\begin{proof}
(i) $\Rightarrow$ (ii). Assume that $\Delta$ has a special cycle $v_1,F_1,...,v_s,F_s,v_{s+1}=v_1$ with $s \ge 3$. Let $\Gamma$ be the 
subcomplex with the facets $F_1,...,F_s$ and $F_1$ a leaf of $\Gamma$. Then there exists a facet $F_i \neq F_1$ such that $F_i \cap F_1 
\neq \emptyset$ and $F_j \cap F_1 \subseteq F_i \cap F_1$ for all $j \neq 1$. Therefore, $v_1,v_2 \in F_i$. Since $F_1$ is the only 
facet of the cycle which contains $v_1,v_2$, we get $F_i = F_1$, a contradiction.\par

(ii) $\Rightarrow$ (iii). Let $M = (a_{ij})$ be a canonical form of the incidence matrix of $\D$.
Assume to the contrary that $M$ has a submatrix of the form $B$.
Since $M$ is canonical, this matrix is contained in a submatrix of $M$ of the form
\[
\begin{pmatrix}
1 & 1 & 0 \\
1 & 0 & 1 \\
0 & 1 & a
\end{pmatrix}.
\]
If $a = 1$, this matrix corresponds to a special cycle of length 3. Since $\D$ has no special cycle of length $\ge 3$, we must have $a 
= 0$. Again, since $M$ is canonical, the new submatrix is contained in a submatrix of $M$ of the form
\[
\begin{pmatrix}
1 & 1 & 0 & 0\\
1 & 0 & 1 & 0\\
0 & 1 & 0 & 1\\
0 & 0 & 1 & b
\end{pmatrix},
\]
we can argue as above to see that $b = 0$. Since we can not go on infinitely, we get a contradiction.\par

(iii) $\Rightarrow$ (i). Since the incidence matrix of every subcomplex of $\D$ is a submatrix of $M$, they are also canonical and have 
no submatrix of the form $B$. Therefore, it suffices to show that $\D$ has a leaf.
Let $M = (a_{ij})$. For each $j = 2,...,m$ let $q_j$ be the smallest integer for which $a_{q_j1} = a_{q_jj} = 1$.
Let $\F(\Delta) = \{F_1, \ldots, F_m\}$, where the order comes from the ordering of the rows of $M$. Since $M$ has no matrix of the 
form $B$, $F_1 \bigcap F_j = F_1 \bigcap \{x_{q_j}, x_{q_j + 1}, \ldots, x_n \}$. Let $j_2, \ldots, j_m$ be a permutation of $2, 
\ldots, s$ such that $q_{j_2} \leq \cdots \leq q_{j_m}$. Then
\[
F_1 \bigcap F_{j_m} \subset \cdots \subset F_1 \bigcap F_{j_2}.
\]
Hence $F_1$ is a leaf of $\Delta$.
\end{proof}

The equivalence of (ii) and (iii) is well-known in hypergraph theory (Hoffman, Kolen and Sakarovitch \cite{Ho}, Lubiw \cite{Lu}) where 
a hypergraph without simple cycle of length $\ge 3$ is called {\it totally balanced} and a matrix of the form in (iii) is called 
greedy. They play an essential role in the study of  Gr\"obner bases of toric ideals arising from bipartite graphs \cite{Hi}. That 
these conditions are equivalent to (i) seems to be not known.

\begin{Corollary} \label{forest standard}
Assume that $\D$ is a forest. Then $A(\D)$ and $A(\D^*)$ are standard graded algebras.
\end{Corollary}

\begin{proof}
This follows from Theorem \ref{balanced} and Theorem \ref{facet ideal 2}
\end{proof}

In \cite[Corollary 3.12]{Fa} Faridi proved that the Rees algebra of the facet ideal $I(\D)$ of a tree $\D$ is a normal Cohen-Macaulay 
domain. By Corollary \ref{facet ideal 1}, this is only a consequence of the fact that $A(\D^*)$ is a standard graded algebra and we 
even obtain stronger properties such as $I(\D)$ is normally torsionfree.

The proof of Theorem \ref{forest} also reveals an interesting relationship between forests and a special kind of leaves.

A leaf $F$ of a simplicial complex $\Delta$ is called a {\em good leaf} if $F$ is a leaf of each subcomplex of $\Delta$ to which it 
belongs. Equivalently, $F$ is a good leaf if the collection of sets $F\sect G$ with $G\in {\mathcal F}(\Delta)$ is totally ordered with 
respect to inclusion. This notion was introduced in the thesis of Zheng \cite{Zh}.

\begin{Corollary} \label{good leaf 1}
\label{goodleaf} Every forest has a good leaf.
\end{Corollary}

\begin{proof}
This follows immediately from  the proof of Theorem \ref{forest} (iii) $\Rightarrow$ (i) ($F_1$ is in fact a good leave).
\end{proof}

If the facets $F_1,\ldots,F_m$ of a simplicial complex $\D$ can be ordered  in such a way that $F_i$ is a good leaf of the subcomplex 
with the facets $F_1,\ldots,F_{i}$ for $i=2,\ldots,m$, we say that $F_1,\ldots,F_m$ is a {\em good leaf order} of $\Delta$. Using 
Corollary \ref{good leaf 1} it follows immediately that a simplicial complex $\Delta$ admits a good leaf order if and only if $\Delta$ 
is  a forest.

It was shown recently by Pelsmajer, Tokaz and D.B. West \cite[Theorem 3.3]{Pe} that totally balanced hypergraphs have good leaf orders, 
where a good leaf is named a simple edge. However, their proof is complicated and needs some backgrounds.

Now we give some nice properties of good leaves.

Given a monomial $m=\prod_{j=1}^rx_{i_j}^{a_i}$. We say this product presentation of $m$ is standard, if
$i_1<i_2<\cdots < i_r$ and $a_1>0,\ldots, a_r>0$.  Of course, if we change the numbering of the variables the standard presentation of 
$m$ also changes. In the following, unless otherwise stated, we always write the monomials
in standard form. A sequence of monomials $m_1,\ldots,m_s$ is said to be an $M$-{\em sequence}\index{M-sequence} if for all $1\leq i 
\leq s$ there exists a numbering of the variables such that if $m_i=x_{i_1}^{a_1}\cdots x_{i_r}^{a_r}$  and whenever $x_{i_k}|m_j$ for 
some  $1\leq k \leq r$ and $i<j$, then $x_{i_k}^{a_k}\cdots x_{i_r}^{a_r}| m_j$. Notice that the numbering of the variables  may depend 
on the index $i$.

It was shown in \cite[Proposition 3.11]{Zh} that if $F_1,\ldots,F_s$ is a good leaf order of a forest, then $x^{F_s},
x^{F_{s-1}},\ldots, x^{F_1}$ is an $M$-sequence. Combining this fact with the result
\cite[Theorem 2.5]{Co} of Conca and De Negri we obtain

\begin{Corollary}
\label{Rees CM} Let $I\subset S$ be a graded ideal. Suppose there exists a term order such that the initial ideal $\ini(I)$ is the 
facet ideal of  a forest. Then
\par {\rm (i) } the Rees algebra of $I$ is Cohen--Macaulay,
\par {\rm (ii) } the associated  graded ring of $I$ is reduced.
\end{Corollary}

Let $d(A(\D))$ denote the maximal degree of the generators of the vertex cover algebra $A(\D)$. The following result shows that 
$d(A(\D))$ does not change by removing a good leaf.

\begin{Theorem} \label{good leaf 2}
Let $\D$ be a simplicial complex and $F$ a good leaf of $\D$. Let $\G$ be the subcomplex obtained from $\D$
by removing $F$. Then
$d(A(\D)) \le d(A(\G))$.
\end{Theorem}

\begin{proof}
By the definition of good leaf, $F$ must contain a vertex which does not belong any facet of $\G$.
Let $1,...,m$ be the vertices of the facets of $\G$. For all $k \ge 1$, if
$(c_1,\cdots,c_m,0,\ldots,0)$ is an indecomposable $k$-cover of $\G$, then
$(c_1,\cdots, c_m,k,\ldots,0)$ is an indecomposable $k$-cover of $\Delta'$. Therefore, $d(A(\D)) \ge d(A(\G))$.

Let $\a=(a_1,\ldots, a_n)$ now be a cover of $\D$ of order $r > d(A(\D))$. We want to find a
vertex cover $\b$ of $\D$ of order $s < r$ such that $\a-\b$ is an $(r-s)$-cover. Then this implies $d(A(\D)) \leq
d(A(\G))$.

Let $\D$ be a simplicial complex on the vertex set $[n]$.
We may assume that $F=\{\ell,\ldots, n\}$ with free vertices
$m,m+1,\ldots, n$ in $\D$. Since $F$ is a good leaf of
$\D$, the set $\{F\sect G:\  G\in\F(\G)\}$ is totally
ordered with respect to inclusion, and we may assume that this set
is the set $\{\{\ell,\ldots, s_1\},\{\ell,\ldots, s_2\},\ldots, \{\ell,\ldots,
s_t\}\}$ with $m-1=s_1>s_2>\cdots >s_t\geq \ell-1$.

If $\ell=m$, then $F\sect G=\emptyset$ for all $G\in\F(\G)$,
and the assertion is trivial. We now may assume that $\ell<m$.

The vector $\c =(a_1,\ldots, a_{m-1},0,\ldots, 0)$ is an $r$-cover of $\G$. Since $r>d(A(\D))$ there exists an
$s$-cover $\e=(e_1,\ldots, e_{m-1},0,\ldots, 0)$ of $\G$ with
$s<r$ such that $\c-\e$ is an $(r-s)$-cover of $\G$.

We are going to construct the desired  $s$-cover $\b$ of
$\D$ by modifying the vertex cover $\e$ of $\G$.
We set $\rho = \sum_{i=\ell}^{m-1}e_i$, and distinguish two cases.

\noindent
Case 1: $\rho<s$. Let $p=\max\{j\in\{\ell,\cdots, n\}:\ \sum_{i=\ell}^{j-1}e_i+\sum_{i=j}^n a_i\geq s\}$. Since
$\sum_{i=\ell}^n a_i\geq  r > s$, such a $j$ exist. We define
$\b=(b_1,\ldots, b_n)$ to be the vector with
\[
b_i=
\begin{cases}
e_i,\quad \text{if}\quad i<p,\\
a_i,\quad \text{if}\quad i>p,\\
\end{cases}
\]
and $b_p=s-\sum_{i = \ell,i\neq p}^n b_i$. Then $\b$ is an $s$-cover of $\D$.
We claim $\a - \b$ is an $(r-s)$-cover of $\D$.
Since $\sum_{i \in F} a_i \geq r$ and since
$\sum_{i \in F} b_i = s$, it is clear that $\sum_{i \in F}
(a_i-b_i) \geq r - s$.  Thus, in particular, since $\sum_{i \in
F} (a_i-b_i) = \sum_{i=\ell}^{p} (a_i - b_i)$, it follows that
$\sum_{i=\ell}^{p} (a_i - b_i) \geq r - s$.  Let $G$ be a facet
of $\G$ and $F \bigcap G = \{\ell,\ldots, s_j \}$. If $s_j <
p$, then $\sum_{i \in G} (a_i-b_i) = \sum_{i \in G} (a_i-e_i)
\geq r - s$.  If $s_j \geq p$, then $\sum_{i \in G} (a_i-b_i)
\geq \sum_{i=\ell}^{p} (a_i-b_i) \geq r - s$, as desired.

\noindent Case 2: $\rho \geq s$. Let
$p=\min\{j\in\{\ell,\cdots, m\}:\ \sum_{i=\ell}^{j}e_i \geq s\}$.
We define $\b=(b_1,\ldots, b_n)$ to be the vector with
\[
b_i=
\begin{cases}
e_i,\quad \text{if}\quad i<p,\\
0,\quad \text{if}\quad i>p, \\
\end{cases}
\]
and $b_p=s-\sum_{i=\ell}^{p-1}b_i$. Then $\a - \b$ is an $(r-s)$-cover of $\D$.
We claim $\b$ is an $s$-cover of $\D$.
One has $\sum_{i \in F} b_i = s$.  Thus, in particular, since
$\sum_{i \in F} b_i = \sum_{i=\ell}^{p} b_i$, it follows that
$\sum_{i=\ell}^{p} b_i = s$.  Let $G$ be a facet of $\G$ and
$F \bigcap G = \{\ell,\ldots, s_j \}$. If $s_j < p$, then $\sum_{i
\in G} b_i = \sum_{i \in G} e_i \geq s$.  Let $s_j \geq p$. Then
$\sum_{i \in G} b_i \geq \sum_{i=\ell}^{p} b_i = s$, as desired.
\end{proof}

\section{Quasi-forests}

A simplicial complex $\Delta$ is called a {\em quasi-forest}, if
the facets of $\Delta$ can be ordered $F_1,\ldots, F_m$ such that
$F_i$ is a leaf of the simplicial complex with facets $
F_1,\ldots, F_{i}$ for $i=2,\ldots,m$. Such an order of the leaves
is called a {\em leaf order.} If in addition, $\Delta$ is
connected then $\Delta$ is called a {\em quasi-tree}.

The name quasi-forest  for this kind of simplicial complexes was
introduced by Zheng \cite{Zh1}. It is clear that any forest (tree)
is a quasi-forest (quasi-tree). The converse is not true, as the
simplicial complex $\Delta$ in Figure \ref{F4} demonstrates.
Indeed, $\Delta$  is a quasi-tree but not a tree.  Because if we
remove the inside facet, then the remaining simplicial complex has
no leaves anymore. In particular, none of the leaves of $\Delta$
is a good leaf.

\begin{figure}[hbt]
\begin{center}
\psset{unit=1cm}
\begin{pspicture}(0,0.5)(3,2.5)
\pspolygon[style=fyp, fillcolor=medium](0.5,0.5)(1.5,0.5)(1,1.4)
 \pspolygon[style=fyp, fillcolor=medium](1.5,0.5)(1,1.4)(2,1.4)
 \pspolygon[style=fyp, fillcolor=medium](1.5,0.5)(2,1.4)(2.5,0.5)
  \pspolygon[style=fyp, fillcolor=medium](1,1.4)(1.5,2.3)(2,1.4)
\end{pspicture}
\end{center}
\caption{}\label{F4}
\end{figure}
\medskip

In hypergraph theory, quasiforests  are  known under the name co-arboreal hypergraphs or hyperforests. We refer to \cite[Chap. 5, \S 
4]{Be2} and \cite{Zha} for more information on this class of hypergraphs.

The significance of quasi-trees results from Dirac's
characterization of chordal graphs. Let $G$ be a finite graph on
$[n]$ without loops, multiple edges and isolated vertices, and
$E(G)$ its edge set.  Dirac proved that $G$ is chordal (i.e.\
every cycle in the graph of length $>3$ has a chord) if and only
if $G$ has a perfect elimination ordering on its vertices. Recall
that a {\em perfect elimination ordering} (or a {\em simplicial
elimination ordering}) is an ordering $v_n,...,v_2,v_1$ on the
vertices of $G$ such that $v_i$ is a simplicial vertex in the
graph induced on vertices $\{v_1,...,v_i\}$.  Here a {\em
simplicial vertex} in a graph is one whose neighbors form a
clique, and clique is a subset $F$ of $[n]$ such that $\{ v, w \}
\in E(G)$ for all $v, w \in F$ with $v \neq w$.

In the terminology of quasi-forests, Dirac's theorem \cite{D} can
be phrased as follows. The following conditions are equivalent:

\begin{enumerate}
\item[(a)] $G$ is a chordal graph. \item[(b)] $G$ is the
$1$-skeleton of a quasi-forest.
\end{enumerate}

In fact, the quasi-forest belonging to the chordal graph $G$ is
the clique complex $\Delta(G)$ of $G$, that is, the simplicial
complex of all cliques of $G$. In \cite{HHZ} there is also given a
higher dimensional version of Dirac's theorem.

A perfect elimination ordering of a graph G is a {\em strong
perfect elimination ordering}  if for all $i<j<k<l$ such that
$\{v_i, v_k\} , \{v_i, v_l\} , \{v_j, v_k\}$ are edges, then
$\{v_j, v_l\}$ is an edge. A graph is {\em strongly chordal} if it
has a strong perfect elimination ordering.

In analogy to Dirac's theorem one has that the following
conditions are equivalent:
\begin{enumerate}
\item[(a)] $G$ is a strongly chordal graph. \item[(b)] $G$ is the
$1$-skeleton of a forest.
\end{enumerate}

This was shown by Farber \cite[Theorem 5.3]{Farber}. He also
showed that strongly chordal graphs can be described by forbidden
induced subgraphs, called trampolines \cite[Theorem 4.1']{Farber}.
For example the 1-skeleton of the simplicial complex in Figure 4
is a trampoline.

While the vertex cover algebra of a forest is standard graded, as
we have seen in Corollary \ref{forest standard}, this is not the
case for quasi-forests. Indeed, the vertex cover of the simplicial
complex in Figure \ref{F4} which assigns to each non-free vertex
the value 1 is an indecomposable vertex cover of order 2. So one
might expect that a quasi-forest which is not a forest has always
indecomposable  vertex covers of higher order. But again this is
not the case. The following figure shows a quasi-tree  which is
not a tree but whose vertex cover algebra is nevertheless standard
graded. Its only non-leaf is the inside 3-dimensional simplex.

\begin{figure}[hbt]
\begin{center}
\psset{unit=0.5cm}
\begin{pspicture}(-1,0.2)(7,6)
\pspolygon[style=fyp,fillcolor=medium](1,3)(0,0)(3,1)
\pspolygon[style=fyp, fillcolor=medium](3,1)(6,0)(5,3)
\pspolygon[style=fyp,fillcolor=medium](5,3)(6,6)(3,5)
\pspolygon[style=fyp, fillcolor=medium](3,5)(0,6)(1,3)
\psline(3,5)(3,1) \psline[linestyle=dotted](1,3)(5,3)
\end{pspicture}
\end{center}
\caption{}\label{Fig6}
\end{figure}

The main goal of this section will be to show that for a certain
restricted class of quasi-forests the vertex cover algebra is
standard graded if and only if the quasi-forest is a forest. For
this purpose we first prove

\begin{Proposition}
\label{need} Let $\Delta$ be a quasi-forest on the vertex set
$[n]$ , and suppose that $\Delta$ contains the subcomplex $\Gamma$
with facets
\begin{align*}
B & =\{i_1,\ldots, i_k\},\;\;  F_a= \{a,i_2, \ldots, i_k\},\\
F_b & =\{i_1,b,i_3, \ldots, i_k\}\;  \text{ and }\; F_c=\{i_1,i_2, c,
i_4,\ldots, i_k\}.
\end{align*}
Then $d(A(\Delta))>1$.
\end{Proposition}

\begin{proof}
Let  $a=(a_1,\ldots,a_n)$ be the integer vector with
\[
a_i=
\begin{cases}
0,\quad \text{if}\quad i\in\{a,b,c,i_4,\ldots, i_k\},\\
1,\quad \text{if}\quad i\in\{i_1,i_2,i_3\},\\
2,\quad \text{otherwise}.
\end{cases}
\]
We show by induction on $\delta=|\Delta|-|\Gamma|$ that
$d(A(\Delta))>1$. If $\delta=0$, then $\Delta=\Gamma$, and $a$ is
obviously an indecomposable  $2$-cover of $\Gamma$.

Now suppose $\delta>0$, and that there is a leaf $F$ of $\Delta$
belonging to $\Delta$ but not to $\Gamma$. Let $\Sigma$ be the
subcomplex of $\Delta$ which is obtained from $\Delta$ by removing
$F$. It is shown in \cite[]{} that $\Sigma$ is again a
quasi-forest, and of course $\Gamma\subset \Sigma$. By induction
hypothesis we have $d(A(\Sigma))>1$. Applying Theorem \ref{good
leaf 2} we conclude that $d(A(\Delta))>1$.

On the other hand, if all leaves of $\Delta$ belong to $\Gamma$ we
will show that $\Delta=\Gamma$, and we obtain again the desired
conclusion.

To see why $\Gamma=\Delta$,  we consider a relation forest
$T(\Delta)$ of $\Delta$ introduced in \cite{HHZ}. The vertices of
$T(\Delta)$ are the facets of $\Delta$. The edges are obtained
recursively as follows: Choose a leaf  $F$ of $\Delta$ and a
branch $G$ of $F$. Then $\{F,G\}$ is an edge of $T(\Delta)$.
Remove $F$ from $\Delta$ and proceed with the remaining
quasi-forest as before to find the other edges of $\Delta$. The
resulting graph may depend on the order of how one chooses the
leaves, but in any case it is a forest. Moreover, each free vertex
of $T(\Delta)$ is a leaf of $\Delta$.

Suppose now that all the leaves of $\Delta$ belong to $\Gamma$.
Since $\Delta$ has at least two leaves  and $B$ is obviously not a
leaf,  we may assume the $F_a$ and $F_b$ are leaves of $\Delta$.
Since $F_a$ and $F_b$ have a face of codimension $1$ in common
with $B$, it follows that $B$ is a branch of $F_a$ and $F_b$ in
$\Delta$. Hence we may construct a relation forest $T(\Delta)$ of
$\Delta$ with edges $\{F_a,B\}$ and $\{F_b,B\}$, and free vertices
$F_a$ and $F_b$.

Since $F_c\not\in\{F_a,F_b,B\}$, the relation forest $T(\Delta)$
must have a free vertex $G\neq F_a,F_b$. If $G\neq F_c$, then
$\Delta$ has a leaf which does not belong to $\Gamma$, a
contradiction. So $G=F_c$, and hence $F_c$ is a leaf  of $\Delta$
and $B$ a branch of $F_c$ in $\Delta$. Thus we add the edge
$\{F_c,B\}$. Therefore $T(\Delta)$ contains the tree $T_0$ given
by the next figure.

\begin{figure}[hbt]
\begin{center}
 \psset{unit=1.0cm}
\begin{pspicture}(0,0.5)(3,2.5)
 \psline(0.52,1.5)(1.43,1.5)
 \psline(1.55,1.45)(2.15,0.85)
 \psline(1.56,1.56)(2.15,2.15)
 \rput(0.45,1.5){$\circ$}
 \rput(1.5,1.5){$\circ$}
 \rput(2.2,0.79){$\circ$}
 \rput(2.2,2.2){$\circ$}
 \rput(0.45,1){$F_a$}
 \rput(1.4,1){$B$}
 \rput(2.5,2.5){$F_b$}
 \rput(2.5,0.5){$F_c$}
\end{pspicture}
\end{center}
\caption{}\label{Fig5}
\end{figure}

Since $T(\Delta)$ has no other free vertices, we must have that
$T(\Delta)=T_0$, and hence $\Delta=\Gamma$.
\end{proof}

Now we have

\begin{Theorem}
\label{also} Let $\Delta$ be a quasi-forest satisfying
\begin{enumerate}
\item[(i)] each connected component of $\Delta$ is connected in
codimension 1; \item[(ii)] each face of codimension 1 belongs to
at most two facets.
\end{enumerate}
Then $\Delta$ is forest if and only if $A(\Delta)$ is standard
graded.
\end{Theorem}

\begin{proof} In view of Corollary \ref{forest standard} it suffices to
show that the vertex cover algebra of a quasi-forest $\Delta$
satisfying the conditions (i) and (ii) is not standard graded,
unless it is a forest. We may assume that $\Delta$ is connected.

Let $F_1,\ldots, F_m$ be a leaf order of  $\Delta$, and set
$\Delta_j=\langle F_1,\ldots, F_j\rangle$ for $j=1,ldots,m$.
Notice that $\Delta_j$  is connected in codimension 1, since this
property is preserved by removing a leaf.

Assuming that $\Delta$ is not a tree, there exists an integer $i$
such that $F_i$ is not a good leaf of $\Delta_i$. In other words,
there exist facets $G$ and $H$ of $\Delta_i$ such that the faces
$G\sect F_{i}$ and $H\sect F_{i}$ are not contained in each other.
In particular, $G\neq H$.  Since $\Delta_{i-1}$ is connected in
codimension 1 there exists a sequence of facets $G=G_0,G_1,\ldots,
G_p=B$ in $\Delta_{i-1}$ where $B$ is the branch of $F_{i}$ such
that $G_k\sect G_{k+1}$ is a face of codimension 1 for
$k=0,\ldots, p-1$. The sequence can be chosen such that $G_s\neq
G_t$ for $s\neq t$.

Similarly there exists a chain $H=H_0,H_1,\ldots, H_q=B$ in
$\Delta_{i-1}$ with $H_s\neq H_t$ if $s\neq t$ and such that
$H_k\sect H_{k+1}$ is a face of codimension 1 for $k=0,\ldots,
q-1$.

 It is clear that $G_{p-1}$ and $H_{q-1}$ are
different from $F_{i}$. We set $G_{p+1}=H_{q+1}=F_{i}$. Then there
exists an integer $r\geq 0$ such that  $G_{p-j}=H_{q-j}$ for all
$0\leq j\leq r$ and $G_{p-r-1}\neq H_{q-r-1}$. Then by assumption
(ii) the simplicial complex $\Gamma$ with facets $G_{p-r+1},
G_{p-r}, G_{p-r-1}$ and $H_{q-r-1}$ is  of the form as described
in Proposition \ref{need},  and hence $A(\Delta)$ is not standard
graded.
\end{proof}


\begin{thebibliography}{1}

\bibitem{Be1} C.\ Berge,   Sur certains hypergraphes généralisant les graphes bipartites.  (French)
Combinatorial theory and its applications, I (Proc. Colloq., Balatonfüred, 1969), pp. 119--133.
North-Holland, Amsterdam, 1970.

\bibitem{Be2} C.\ Berge, Hypergraphs. Combinatorics of finite sets. North-Holland, Amsterdam, 1989.

\bibitem{Be3} C.\ Berge and M.\ Las Vergnas, Sur un théorème du type König pour hypergraphes.  Ann. New York Acad. Sci.\  {\bf 175} 
(1970), 32--40.

\bibitem{Br} W. Bruns and J. Herzog, Cohen-Macaulay rings. Cambridge University Press, Cambridge, 1993.

\bibitem{Co} A.\ Conca and M.\ DeNegri,$M$-sequences, graph ideals, and ladder ideals of linear type.  J. Algebra  {\bf 211}  (1999),  
no. 2, 599--624.

\bibitem{D} G.\ A.\ Dirac, On rigid circuit graphs, Abh.\ Math.\ Sem.\ Univ.\ Hamburg, {\bf 38} (1961), 71--76.

\bibitem{Du} H.\ Duchet, Hypergraphs.  Handbook of combinatorics, Vol. 1, 381--432, Elsevier, Amsterdam, 1995.

\bibitem{Es} C.A. Escobar, R. Villarreal and Y. Yoshino, Torsion freeness and normality of blowup rings of monomial ideals.  
Commutative algebra,  69--84, Lect. Notes Pure Appl. Math. 244, Chapman \& Hall, 2006.

\bibitem{Farber} M.\ Farber, Characterizations of strongly chordal graphs,  Discrete Math.\  {\bf  43}  (1983), no. 2-3, 173--189.

\bibitem{Fa} S. Faridi, The facet ideal of a simplicial complex,  Manuscripta Math. {\bf 109} (2002), 159 -- 174.

\bibitem{Fu} D.R. Fulkerson, A.J. Hoffman, R. Oppenheim, On balanced matrices.  Math. Programming Stud.\ {\bf 1} (1974), 120--132.

\bibitem{Gi1} I. Gitler, E. Reyes and R. Villarreal, Blowup algebras of ideals of vertex covers of bipartite graphs.  Algebraic 
structures and their representations,  273--279, Contemp. Math., 376, Amer. Math. Soc., Providence, RI, 2005.

\bibitem{Gi} I. Gitler, C.E. Valencia and R. Villarreal, A note on Rees algebras and the MFMC property, Preprint 2005.

\bibitem{He} J.\ Herzog, T.\ Hibi and N.V.\ Trung, Symbolic powers of monomial ideals and vertex cover
algebras, to appear in Adv. Math. (ArXiv, math.AC/0512423).

\bibitem{HHZ} J.\ Herzog, T.\  Hibi and X.\ Zheng, Dirac's theorem on chordal
graphs and Alexander duality. European J.\ Comb. 25(7) (2004),
949--960.

\bibitem{Hi}  T.\ Hibi and H. Ohsugi,  Koszul bipartite graphs, 
Advances in Applied Math. {\bf 22} (1999), 25 -- 28.

\bibitem{Hoch} M.\ Hochster, Rings of invariants of tori, Cohen--Macaulay rings generated by
monomials, and polytopes. Ann. Math. {\bf 96} (1972), 318--337.

\bibitem{Ho} A.J. Hoffman, A.W.J. Kolen and M. Sakarovitch,  Totally-balanced and greedy matrices.
SIAM J. Algebraic Discrete Methods 6 (1985), no. 4, 721--730.

\bibitem{Hu} C. Huneke, On the associated graded ring of an ideal.  Illinois J. Math.  {\bf 26}  (1982), no. 1, 121--137.

\bibitem{HuV} C. Huneke, B. Ulrich and V. Vasconcelos, Reduced normal cones are domains.  Invariant theory (Denton, TX, 1986),  
95--101, Contemp. Math. 88, Amer. Math. Soc., Providence, RI, 1989.

\bibitem{Lu} A. Lubiw, Doubly lexical orderings of matrices.  SIAM J. Comput.  16  (1987),  no. 5, 854--879.

\bibitem{Pe} M. Pelsmajer, J. Tokaz and D.B. West, New proofs for strongly chordal graphs and chordal bipartite graphs, Preprint, 2004.

\bibitem{Si} A. Simis, V. Vasconcelos and Villarreal, On the ideal theory of graphs.  J. Algebra  {\bf 167}  (1994),  no. 2, 389--416.

\bibitem{Sta} R. Stanley, Combinatorics and Commutative Algebra (second edition), Birkh\"auser, 1996.

\bibitem{Zha}  L.W. Zhang, Studies on hypergraphs. I. Hyperforests. Discrete Appl. Math. {\bf 42} (1993), no. 1, 95--112.

\bibitem{Zh} X.\ Zheng, Homological properties of monomial ideals associated to
quasi-trees and lattices, Dissertation, Essen,  August 25, 2004.

\bibitem{Zh1} X.\ Zheng, Resolutions of facet ideals. Comm. Algebra {\bf 32}, no. 6, (2004),
2301--2324.
\end{thebibliography}
\end{document}